\documentclass[twoside]{article}
\usepackage{graphicx}
\usepackage{amsmath, amsthm, amssymb,mathrsfs}

\def \hs {\hspace*{-2mm}}
\newfont{\mi}{cmti9}
\newfont{\m}{cmr8}
\newfont{\ms}{cmsl8}
\newfont{\autor}{cmcsc10}
\newtheorem{theorem}{Theorem}
\newtheorem{remark}{Remark}
\newtheorem{lemma}{Lemma}
\newtheorem{corollary}{Corollary}
\begin{document}
\rule[3mm]{128mm}{0mm} \vspace*{-16mm}

{\footnotesize
K\,Y\,B\,E\,R\,N\,E\,T\,I\,K\,A\, ---
\,V\,O\,L\,U\,M\,E\, {\it X\,X\,}
(\,2\,0\,0\,*\,)\,,\,
N\,U\,M\,B\,E\,R\, X\,,\,\
P\,A\,G\,E\,S\, \,x\,x\,x -- x\,x\,x}\\
\rule[3mm]{128mm}{0.2mm}

\vspace*{11mm}

{\large\bf \noindent Optimal\,
Sequential\, Multiple\, Hypothesis\,
Testing\, in\, Presence\, of\,
Control\, Variables}

\vspace*{8mm}

{\autor \indent Andrey Novikov}\\

\vspace*{23mm}

\small

Suppose that at any stage of a
statistical experiment  a  control
variable $X$ that affects the
distribution of the observed data $Y$
at this stage can be used. The
distribution of $Y$ depends on some
unknown parameter $\theta$, and  we
consider the problem of testing
multiple hypotheses
$H_1:\,\theta=\theta_1$,
$H_2:\,\theta=\theta_2, \dots$,
$H_k:\,\theta=\theta_k$ allowing the
data to be controlled by $X$, in the
following sequential context.
 The experiment starts with assigning a value $X_1$ to the control
variable and observing $Y_1$ as a
response. After some analysis, another
value $X_2$ for the control variable is
chosen, and  $Y_2$ as a response is
observed, etc.  It is supposed that the
experiment eventually stops, and at
that moment a final decision in favor
of one of the hypotheses $H_1,\dots$,
$H_k$ is to be taken. In this article,
our aim  is to characterize the
structure of optimal sequential testing
procedures based on data obtained from
an experiment of this type in the case
when the observations $Y_1, Y_2,\dots,
Y_n$ are independent, given controls
$X_1,X_2,\dots, X_n$,  $n=1,2,\dots$.

\smallskip\par
\noindent {\sl Keywords:}\,
\begin{minipage}[t]{112mm}
sequential analysis, sequential
hypothesis testing, multiple
hypotheses, control variable,
independent observations,  optimal
stopping, optimal control, optimal
decision, optimal sequential testing
procedure, Bayes, sequential
probability ratio test
\end{minipage} \smallskip
\par
\noindent {\sl AMS Subject
Classification:} 62L10,  62L15, 60G40,
62C99, 93E20

\normalsize


\section{\hs  INTRODUCTION. PROBLEM SET-UP.}\label{s1}
Let us suppose that at any stage of a statistical experiment  a
"control variable" $X$ can be used, that affects the distribution of
the observed data $Y$ at this stage. "Statistical" means that  the
distribution of $Y$ depends on some unknown parameter $\theta$, and
we
 have the usual goal of statistical analysis: to obtain some information
about the true value of $\theta$. In this work, we consider the
 problem of testing multiple hypotheses
$H_1:\,\theta=\theta_1$, $H_2:\,\theta=\theta_2, \dots$,
$H_k:\,\theta=\theta_k$  allowing the data to be controlled by $X$,
in the following "sequential" context.

 The experiment starts with assigning a value $X_1$ to the control
variable and observing $Y_1$ as a response. After some analysis, we
choose another value $X_2$ for the control variable, and observe
$Y_2$ as a response. Analyzing this, we choose $X_3$ for the third
stage, get $Y_3$, and so on. In this way, we obtain a sequence
$X_1,\dots, X_n$, $Y_1,\dots, Y_n$ of experimental data,
$n=1,2,\dots$. It is supposed that the experiment eventually stops,
and at that moment a final decision in favor of one of $H_1,\dots$,
 $H_k$ is to be taken.

In this article, our aim  is to characterize the structure of
optimal sequential procedures, based on this type of data, for
testing the multiple hypotheses $H_1,\dots$,
 $H_k$.

  We follow \cite{Malyutov1983} and \cite{Volodin}  in our interpretation of "control variables". For example, in a regression experiment, with a dependent variable $Y$ and an independent variable $X$, the variable $X$ is a control variable in our sense, whenever the experimenter can vary its value before the next observation is taken.
  Another classical context for "control variables" in our sense is the experimental design, when one of some alternative treatments is assigned to every experimental unit before the experiment starts. The randomization, which is frequently used with both these type of "controlled" experiments, can be easily incorporated in our theory below as well.

  There exist yet another concept of  "control variables" introduced by Haggstrom \cite{Haggstrom}, and largely used in \cite{LectureNotes} and many subsequent articles (see also \cite{Cressie} for  results, closely related  to \cite{LectureNotes}, where "control variables" are not used). In the context of \cite{LectureNotes}, a control variable, roughly speaking, is an integer variable whose value, at every stage of the experiment, is a prescription of a number of the additional observations to be taken at the next stage, if any. To some extent, it is related to our control variables as well, because it affects the distribution of subsequently observed data. It is very likely that our method will work for this type of "sequentially planned" experiments as well, but formally it does not fit our theory below, mainly because we do not allow that the cost of observations depend on $X$.

 In this article, we follow very closely our article \cite{NovikovIMF}, where the case of $k=2$
 simple hypotheses was considered, and use a method  based on the same ideas as in \cite{NovikovKyb},
 where multiple hypothesis testing for experiments without control variables was studied.

For  data vectors, let us write, briefly, $X^{(n)}$ instead of $(X_1,\dots,X_n)$,
$Y^{(n)}$ instead of $(Y_1,\dots,Y_n)$, $n=1,2,\dots$, etc.
 Let us define a (randomized) sequential hypothesis testing
 procedure
 as a triplet $( \chi,\psi, \phi)$ of a a {\em control
policy} $\chi$, a {\em stopping rule} $\psi$, and a {\em decision
rule} $\phi$,
 with $$
\chi=\left(\chi_1, \chi_2,
\dots,\chi_n,\dots\right),
\quad\psi=\left(\psi_1, \psi_2, \dots
,\psi_{n},\dots\right),\quad
\phi=\left(\phi_1, \phi_2, \dots
,\phi_n,\dots\right),$$ with the
components described below.

The functions
$$\chi_{n}=\chi_{n}(x^{(n-1)},y^{(n-1)}),\quad n=1,2,\dots
$$ are supposed to be measurable functions with values
in the space of values of the control variable.  The  functions
$$\psi_n=\psi_n(x^{(n)},y^{(n)}), \quad n=1,2,\dots$$ are supposed to be some measurable functions with values in $[0,
1]$. Finally,
$$\phi_n=(\phi_{n1},\phi_{n2},\dots, \phi_{nk}),$$
with
$$
\phi_{ni}=\phi_{ni}(x^{(n)},y^{(n)}),\quad i=1,\dots, k,
$$
are supposed to be measurable non-negative functions such that
$$\sum_{i=1}^k\phi_{ni}(x^{(n)},y^{(n)})\equiv 1\quad\mbox{for any}\;n=1,2\dots.$$

The interpretation of all these functions is as follows.

The experiments starts at stage $n=1$ applying $\chi_1$ to determine
the initial control $x_1$. Using this control, the first data $y_1$
is observed.

At any stage $n\geq 1$:  the value of $\psi_n(x^{(n)},y^{(n)})$ is
interpreted as the conditional probability {\em to stop and proceed
to decision making}, given that that we came to that stage and that
the observations  were $(y_1, y_2, \dots, y_n)$ after the respective
controls $(x_1,x_2,\dots,x_n)$ have been applied. If there is no
stop, the experiments continues to the next stage ($n+1$), defining first
the new control value $x_{n+1}$ by applying the control policy:
$$
x_{n+1}=\chi_{n+1}(x_1, \dots, x_{n};y_1,\dots, y_{n})
$$
and then taking an additional observation $y_{n+1}$ using control
$x_{n+1}$.
Then the rule $\psi_{n+1}$ is applied to $(x_1,\dots,
x_{n+1};y_1,\dots,y_{n+1})$ in the same way as as above, etc., until
the experiment eventually stops.

It is supposed that when the experiment stops,  a decision {\em to
accept one and only one} of $H_1,\dots$, $H_k$ is to be made. The
function $\phi_{ni}(x^{(n)},y^{(n)})$ is interpreted  as the
conditional probability {\em to accept} $H_i$, given that the
experiment stops at stage $n$ being  $(y_1,\dots, y_n)$ the data
vector observed and $(x_1,\dots,x_n)$ the respective controls
applied.

The control policy $\chi$ generates, by the above process, a
sequence of random variables $X_1,X_2,\dots,X_n$, recursively by
$$
X_{n+1}=\chi_{n+1}(X^{(n)},Y^{(n)}).
$$
 The stopping rule $\psi$ generates, by the above
process, a random variable $\tau_\psi$ ({\em stopping time}) whose
distribution is given by
\begin{equation}\label{2a}
P_\theta^\chi(\tau_\psi=n)=E_\theta^\chi (1-\psi_1)(1-\psi_2)\dots
(1-\psi_{n-1})\psi_n.
\end{equation}
Here, and throughout the paper, we interchangeably use  $\psi_n$
both for $$\psi_n(x^{(n)},y^{(n)})$$ and for
$$\psi_n(X^{(n)},Y^{(n)}),$$ and so do we for any other
function  $$F_n=F_n(x^{(n)},y^{(n)}).$$ This does not cause any
problem if we adopt the following agreement: when $F_n$ is under
probability or expectation sign, it is $F_n(X^{(n)},Y^{(n)})$,
otherwise it is $F_n(x^{(n)},y^{(n)})$.


For a sequential testing procedure $(\chi,\psi,\phi)$ let us define
\begin{equation}\label{3a}
\alpha_{ij}(\chi,\psi,\phi)=P_{\theta_i}(\,\mbox{accept}\,
H_j)=\sum_{n=1}^\infty E_{\theta_i}^\chi(1-\psi_1)
\dots(1-\psi_{n-1})\psi_n\phi_{nj}
\end{equation}
and
\begin{equation}\label{3c}
\beta_{i}(\chi,\psi,\phi)=P_{\theta_i}(\,\mbox{accept any}\,
H_j\,\mbox{different from}\,
H_i)=\sum_{j\not=i}\alpha_{ij}(\chi,\psi,\phi).
\end{equation}

The probabilities $\alpha_{ij}(\chi,\psi,\phi)$ for $j\not=i$ can be
considered "individual" error probabilities and
$\beta_i(\chi,\psi,\phi)$ "gross" error probability, under
hypothesis $H_i$, $i=1,2,\dots, k$, of the sequential testing procedure
$(\chi,\psi,\phi)$.


Another important characteristic of a sequential testing procedure
is the {\em average sample number}:
\begin{equation}\label{8aa}
N(\theta;\chi,
\psi)=E_\theta^\chi\tau_\psi=\begin{cases}\sum_{n=1}^\infty
nP_\theta^\chi(\tau_\psi=n), \;\mbox{if}\;
P_\theta^\chi(\tau_\psi<\infty)=1,\cr
\infty\quad\mbox{otherwise}.\end{cases}
\end{equation}

 In this article, we solve the two following  problems:
\begin{description}
\item[Problem I.]
Minimize $N(\chi,\psi)=N(\theta_1;\chi,\psi)$  over all sequential
testing procedures $(\chi,\psi,\phi)$ subject to
\begin{equation}\label{1}
\alpha_{ij}(\chi,\psi,\phi)\leq\alpha_{ij},\quad\mbox{for any}\;
i=1,\dots k,\;\mbox{and for any}\; j\not= i,
\end{equation}
where $\alpha_{ij}\in(0,1)$ (with $i,j=1,\dots k$, $j\not=i$) are
some constants.
\item[Problem II.]
Minimize $N(\chi,\psi)=N(\theta_1;\chi,\psi)$  over all sequential testing procedures
$(\chi,\psi,\phi)$ subject to
\begin{equation}\label{2}
\beta_{i}(\chi,\psi,\phi)\leq\beta_i,\quad\mbox{for any}\; i=1,\dots
k,
\end{equation}
with some constants $\beta_i\in(0,1)$, $i=1,\dots,k$.
\end{description}

In Section \ref{s2}, we reduce the problem of minimizing
$N(\chi,\psi)$ under constraints (\ref{1}) (or (\ref{2})) to an
unconstrained minimization problem. The new objective function is
the Lagrange-multiplier function $L(\chi,\psi,\phi)$.

 Then, finding
$$
L(\psi,\phi)=\inf_{\phi}L(\chi,\psi,\phi)
$$
we reduce the problem further  to a problem of finding optimal control policy and stopping rule.

In Section \ref{s3}, we solve the problem of minimizition of  $L(\chi,\psi)$ in a class of
control-and-stopping strategies.

In Section \ref{s4}, the likelihood ratio structure for optimal
strategy is given.

In Section \ref{s5}, we apply the results obtained in Sections
\ref{s2} -- \ref{s4} to the solution of Problems I and II.

The final Section \ref{s6} contains some additional results, examples and discussion.

\section{\hs REDUCTION TO A PROBLEM OF OPTIMAL CONTROL AND STOPPING}\label{s2}

In this section, Problems I and II will be reduced to
unconstrained optimization problems using the idea of the Lagrange
multipliers method.

\subsection{\normalsize Reduction to Non-Constrained Minimization in Problems I and II}
 The following two theorems are practically Theorem 1 and Theorem 2 in \cite{NovikovKyb}. They reduce Problem I and Problem II to respective unconstrained minimization problems, using the idea of the Lagrage multipliers method.

For Problem I, let us define $L(\chi,\psi,\phi)$ as
\begin{equation}\label{4}
L(\chi,\psi,\phi)=N(\chi,\psi)+\sum_{{1\leq i,j\leq
k};\,{i\not=j}}\lambda_{ij}\alpha_{ij}(\chi,\psi,\phi)
\end{equation}
where  $\lambda_{ij}\geq 0$   are some constant multipliers.

Let $\Delta$ be a class of sequential testing procedures.

\begin{theorem}\label{t1} Let exist $\lambda_{ij}> 0$, $i=1,\dots, k$, $j=1,\dots, k$, $j\not=i$, and a testing procedure $(\chi^*,\psi^*,\phi^*)\in
\Delta$ such that for any other testing procedure
$(\chi,\psi,\phi)\in \Delta$
\begin{equation}\label{5}
L(\chi^*,\psi^*,\phi^*)\leq L(\chi,\psi,\phi)
 \end{equation}
 holds (with $L(\chi,\psi,\phi)$ defined by (\ref{4})), and such that
 \begin{equation}\label{6}\alpha_{ij}(\chi^*,\psi^*,\phi^*)=\alpha_{ij}\quad\mbox{for any}\;
i=1,\dots k,\;\mbox{and for any}\; j\not= i.
 \end{equation}

 Then for any testing procedure $(\chi,\psi,\phi)\in\Delta$ such
 that
 \begin{equation}\label{5bis}
 \alpha_{ij}(\chi,\psi,\delta)\leq\alpha_{ij}\quad\mbox{for any}\;
i=1,\dots k,\;\mbox{and for any}\; j\not= i,
 \end{equation}
 it holds
\begin{equation}\label{5a}
N(\chi^*,\psi^*)\leq  N(\chi,\psi). \end{equation}

 The inequality
in (\ref{5a}) is strict if at least one of the equalities
(\ref{5bis}) is strict.
\end{theorem}

For Problem II, let now $L(\chi,\psi,\phi)$ be defined as
\begin{equation}\label{4aBis}
L(\chi,\psi,\phi)=N(\chi,\psi)+\sum_{i=1}^k\lambda_{i}\beta_{i}(\chi,\psi,\phi),
\end{equation}
where $\lambda_i\geq 0$ are the Lagrange multipliers.

\begin{theorem}\label{t1a} Let exist $\lambda_{i}> 0$, $i=1,\dots, k$,  and a testing procedure $(\chi^*,\psi^*,\phi^*)\in
\Delta$ such that for any other testing procedure
$(\chi,\psi,\phi)\in \Delta$
\begin{equation}\label{5Bis}
L(\chi^*,\psi^*,\phi^*)\leq L(\chi,\psi,\phi)
 \end{equation}
 holds  (with $L(\chi,\psi,\phi)$ defined by (\ref{4aBis})), and such that
 \begin{equation}\label{6BisBis}\beta_{i}(\chi^*,\psi^*,\phi^*)=\beta_{i}\quad\mbox{for any}\;
i=1,\dots k.
 \end{equation}

 Then for any testing procedure $(\chi,\psi,\phi)\in\Delta$ such
 that
 \begin{equation}\label{5bisBis}
 \beta_{i}(\chi,\psi,\delta)\leq\beta_{i}\quad\mbox{for any}\;
i=1,\dots k,
 \end{equation}
 it holds
\begin{equation}\label{5aBis}
N(\chi^*,\psi^*)\leq  N(\chi,\psi). \end{equation}

 The inequality
in (\ref{5aBis}) is strict if at least one of the equalities
(\ref{5bisBis}) is strict.
\end{theorem}

\subsection{\hs Optimal Decision Rules}

Due to Theorems \ref{t1} and \ref{t1a}, Problem I is reduced to
minimizing (\ref{4}) and Problem II is reduced to minimizing
(\ref{4aBis}). But (\ref{4aBis}) is a particular case of (\ref{4}),
namely, when $\lambda_{ij}=\lambda_i$ for any $j=1,\dots,k$,
$j\not=i$ (see (\ref{3a}) and (\ref{3c})). Because of that, we will only solve the problem of minimizing $L(\chi,\psi,\phi)$
defined by (\ref{4}).

In particular, in this section  we  find
$$
\inf_{\phi}L(\chi,\psi,\phi),
$$
and the corresponding decision rule $\phi$, at which this infimum is
attained.

 Let $I_A$ be the
indicator function of the event $A$.

From this time on, we suppose that for any $n=1,2,\dots,$ the random
variable $Y$, when a control $x$ is applied, has a probability
"density" function
\begin{equation}\label{0} f_\theta(y|x)
 \end{equation}
 (Radon-Nicodym derivative of its distribution) with respect to
 a $\sigma$-finite measure $\mu$ on the respective space.
 We are supposing as well that, at any stage $n\geq 1$, given control values $x_1,x_2,\dots
 x_n$ applied, the observations $Y_1,Y_2,\dots, Y_n$ are
 independent, i.e. their joint probability density function, conditionally on given controls $x_1,x_2,\dots
 x_n$,  can be calculated as
\begin{equation}\label{32}
 f_\theta^n(x_1,\dots,x_n;y_1,\dots, y_n)=\prod_{i=1}^n f_\theta (y_i|x_i),
 \end{equation}
 with respect to the product-measure
 $\mu^n=\mu\otimes\dots\otimes\mu$ of $\mu$ $n$ times by itself.
It is easy to see that any expectation, which uses a control policy
$\chi$, can be expressed as
$$
E_{\theta}^\chi g(Y^{(n)})=\int
g(y^{(n)})f_\theta^{n,\chi}(y^{(n)})d\mu^n(y^{(n)}),
$$
where
$$
f_\theta^{n,\chi}(y^{(n)})=\prod_{i=1}^n f_\theta (y_i|x_i)
$$
with
\begin{equation}\label{31}x_i=\chi_i(x^{(i-1)},y^{(i-1)})\end{equation} for
any $i=1,2,\dots$.

Similarly, for any function $F_n=F_n(x^{(n)},y^{(n)})$ let us define
$$F_n^\chi(y^{(n)})=F_n(x^{(n)},y^{(n)})$$ where $x_1,\dots, x_n$
are defined by (\ref{31}).

 As a first step of minimization of $L(\chi,\psi,\phi)$, let us prove the following
\begin{theorem} \label{t2} For any $\lambda_{ij}\geq 0$, $i=1,\dots,k$, $j\not= i$,  and for
any sequential testing procedure $(\chi,\psi,\phi)$
\begin{equation}\label{7}
L(\chi,\psi,\phi)\geq N(\chi,\psi)+\sum_{n=1}^\infty \int
(1-\psi_1^\chi)\dots (1-\psi_{n-1}^\chi)\psi_n^\chi l_n^\chi d\mu^n,
\end{equation}
where
\begin{equation}\label{6a}
l_n=\min_{1\leq j\leq k}\sum_{i\not =j}\lambda_{ij}f_{\theta_i}^n.
\end{equation}
The right-hand side of (\ref{7}) is attained if
\begin{equation}\label{7aa}
\phi_{nj}\leq I_{\left\{\sum_{i\not
=j}\lambda_{ij}f_{\theta_i}^n=l_n\right\}}
\end{equation}
for any $n=1,2,\dots$ and for any $j=1,\dots k$.
\end{theorem}
\begin{proof} Let us suppose that $N(\chi,\psi)<\infty$, otherwise (\ref{7}) is trivial.
Then let us prove an equivalent to (\ref{7}) inequality:
\begin{equation}\label{9}
\sum_{1\leq i,j\leq
k;\,j\not=i}\lambda_{ij}\alpha_{ij}(\chi,\psi,\phi)\geq\sum_{n=1}^\infty
\int (1-\psi_1^\chi)\dots (1-\psi_{n-1}^\chi)\psi_n^\chi l_n^\chi
d\mu^n.
\end{equation}

The left-hand
side of it
 can be represented as
\begin{equation}\label{7f}
\sum_{1\leq i,j\leq
k;\,j\not=i}\lambda_{ij}\alpha_{ij}(\chi,\psi,\phi)= \sum_{n=1}^\infty \int (1-\psi_1^\chi)\dots (1-\psi_{n-1}^\chi)\psi_n^\chi
 \sum_{j=1}^k\left(\sum_{1\leq i\leq
k;\,i\not=j}\lambda_{ij}f_{\theta_i}^{n,\chi}\right)\phi_{nj}^\chi
d\mu^n
 \end{equation}
(see (\ref{3a})).

Applying Lemma 1 \cite{NovikovKyb} to each summand on the right-hand side of
(\ref{7f}) we immediately have:
 \begin{equation}\label{7g}
 \sum_{1\leq i,j\leq
k;\,j\not=i}\lambda_{ij}\alpha_{ij}(\chi,\psi,\phi) \geq \sum_{n=1}^\infty \int (1-\psi_1^\chi)\dots (1-\psi_{n-1}^\chi)\psi_n^\chi
 l_n^\chi
d\mu^n
 \end{equation}
with an equality if
$$
\phi_{nj}\leq I_{\{\sum_{i\not=j}\lambda_{ij}f_{\theta_i}^{n}=l_n\}}
$$
for any $n=1,2,\dots$ and for any $1\leq j\leq k$.
\end{proof}

\begin{remark}\label{r2}  It is easy to see, using (\ref{8aa}) and (\ref{7g}), that
\begin{equation}\label{10}
L(\chi,\psi)=\inf_{\phi}L(\chi,\psi,\phi)=\sum_{n=1}^\infty \int
(1-\psi_1^\chi)\dots
(1-\psi_{n-1}^\chi)\psi_n^\chi\left(nf_{\theta_1}^{n,\chi}+l_n^\chi
\right) d\mu^n,
\end{equation}
with $l_n$  defined by (\ref{6a}), if $P_{\theta_1}^\chi(\tau_\psi<\infty)=1$, and $L(\chi,\psi)=\infty$
otherwise.
\end{remark}

Problem I is reduced now to the problem of finding strategies $(\chi,\psi)$ which minimize $L(\chi,\psi)$.
Indeed, if there is a $(\chi^*,\psi^*)$ such that
$$
L(\chi^*,\psi^*)=\inf_{(\chi,\psi)}L(\chi,\psi),
$$
then  for any $\phi^*$ satisfying
$$
\phi_{nj}^*\leq I_{\left\{\sum_{i\not
=j}\lambda_{ij}f_{\theta_i}^n=l_n\right\}}
$$
(see (\ref{7aa})), by Theorem \ref{t2} for any $(\chi,\psi,\phi)$
$$
L(\chi^*,\psi^*,\phi^*)=L(\chi^*,\psi^*)\leq L(\chi,\psi)=L(\chi,\psi,\phi^*),
$$
thus, the conditions of Theorem \ref{t1} are fulfilled with $\alpha_{ij}=\alpha_{ij}(\chi^*,\psi^*,\phi^*)$ for $i,j=1,\dots, k$, $i\not = j$.

Because of this, in what follows we solve the problem of minimizing $L(\chi,\psi)$.

Let us denote, for the rest of this article,
$$
s_n^{\psi}=(1-\psi_1)\dots(1-\psi_{n-1})\psi_n\quad\mbox{and}\quad c_n^{\psi}=(1-\psi_1)\dots(1-\psi_{n-1})
$$
for any $n=1,2,\dots$ (being $s_1^\psi\equiv \psi_1$ and $c_1^\psi\equiv 1$). Respectively,
$$
s_n^{\psi,\chi}=(1-\psi_1^\chi)\dots(1-\psi_{n-1}^\chi)\psi_n^\chi\quad\mbox{and}\quad c_n^{\psi,\chi}=(1-\psi_1^\chi)\dots(1-\psi_{n-1}^\chi)
$$
for any $n=1,2,\dots$ (being $s_1^{\psi,\chi}\equiv \psi_1^\chi$ and $c_1^{\psi,\chi}\equiv 1$ as well).

Let also
$$
C_{n}^{\psi,\chi}=\{y^{(n)}:(1-\psi_1^\chi(y^{(1)}))\dots(1-\psi_{n-1}^\chi(y^{(n-1)}))>0\},
$$
for any $n\geq 2$, and  let $C_1^{\psi,\chi}$ be the space of all  $y^{(1)}$, and finally let
$$
\bar C_{n}^{\psi,\chi}=\{y^{(n)}:(1-\psi_1^\chi(y^{(1)}))\dots(1-\psi_{n}^\chi(y^{(n)}))>0\},
$$
for any $n\geq 1$.

\section{\hs OPTIMAL CONTROL AND STOPPING}\label{s3}

In this section, the problem of finding strategies $(\chi,\psi)$ minimizing
$L(\chi,\psi)$ (see (\ref{10})) will be solved.

\subsection{\hs Truncated Stopping rules}

In this section, we solve, as an intermediate step,  the problem of
minimization of $L(\chi,\psi)$ over all strategies $(\chi,\psi)$ with truncated stopping rules, i.e.
such $\psi$ that
\begin{equation}\label{11}
\psi=(\psi_1,\psi_2,\dots,\psi_{N-1},1,\dots).\end{equation}
Let $\Delta^N$  be the class of stopping rules $\psi$ of type (\ref{11}), where $N$ is any integer, $N\geq 2$.

The following Theorem can be proved in the same way as Theorem 4.2 in \cite{NovikovIMF}.
\begin{theorem}\label{t3}  Let $\psi\in\Delta^N$ be any (truncated) stopping rule, and $\chi$ any control policy. Then for any
$1\leq r\leq N-1$ the following inequalities hold true
\begin{equation}\label{46a}
L(\chi,\psi)\geq\sum_{n=1}^{r}\int
s_n^{\psi,\chi}(nf_{\theta_1}^{n,\chi}+l_n^\chi)d\mu^n
+\int
c_{r+1}^{\psi,\chi}\left((r+1)f_{\theta_1}^{r+1,\chi}+V_{r+1}^{N,\chi}\right)d\mu^{r+1}
\end{equation}
\begin{equation}\label{46b}
\geq \sum_{n=1}^{r-1}\int
s_n^{\psi,\chi}(nf_{\theta_1}^{n,\chi}+l_n^\chi)d\mu^n
+\int
c_r^{\psi,\chi}\left(rf_{\theta_1}^{r,\chi}+V_{r}^{N,\chi}\right)d\mu^r,
\end{equation}
where $V_N^N\equiv l_N$, and recursively for $n=N, N-1, \dots 2 $
\begin{equation}\label{48}
V_{n-1}^N=\min\{l_{n-1},f_{\theta_1}^{n-1}+R_{n-1}^N\},
\end{equation}
with
\begin{equation}\label{48a}
R_{n-1}^N=R_{n-1}^N(x^{(n-1)};y^{(n-1)})=\min_{x_{n}}\int V_{n}^N(x_1,\dots,x_n;y_1,\dots,y_n)d\mu(y_{n}).
\end{equation}

 The lower bound in (\ref{46b}) is attained if and only if
\begin{equation}\label{49}
I_{\{l_{n}^\chi< f_{\theta_1}^{n,\chi}+R_n^{N,\chi} \}}\leq\psi_{n}^\chi\leq I_{\{l_{n}^\chi\leq f_{\theta_1}^{n,\chi}+R_n^{N,\chi} \}}
\end{equation}
$\mu^n$-almost everywhere on $C_n^{\psi,\chi}$
and
\begin{equation}\label{49a}
R_n^{N,\chi}(y^{(n)})=\int V_{n+1}^{N,\chi}d\mu(y_{n+1})
\end{equation}
 $\mu^n$-almost everywhere on $\bar C_{n}^{\psi,\chi}$, for any $n=r,\dots, N-1$.
\end{theorem}
\begin{remark}\label{r4}
It is supposed in Theorem \ref{t3}, and in what follows in this article, that all the functions $R_{n-1}^N$
defined by (\ref{48a}) are well-defined and
measurable, for any $n=2,\dots, N$, and for any $N=1,2,\dots$.
\end{remark}
The following Corollary characterizes optimal strategies with truncated stopping rules. It immediately follows from Theorem \ref{t3} applied for  $r=1$.
\begin{corollary}\label{c1} For any truncated stopping rule $\psi\in\Delta^N$, and for any control rule $\chi$
\begin{equation}\label{49b}
L(\chi,\psi)\geq  1+R_0^N,
\end{equation}
where
\begin{equation}\label{50b}
R_{0}^N=\min_{x_{1}}\int V_{1}^N(x_1;y_1)d\mu(y_{1}).
\end{equation}
The lower bound in (\ref{49b}) is attained if and only if (\ref{49}) is satisfied $\mu^n$-almost everywhere on $C_n^{\psi,\chi}$ and
(\ref{49a}) is satisfied  $\mu^n$-almost everywhere on $\bar C_{n}^{\psi,\chi}$, for any $n=1,2,\dots, N-1$ and,
additionally,
\begin{equation}\label{49e}
R_0^N=\int V_{1}^{N}(\chi_1;y_1)d\mu(y_{1}).
\end{equation}
\end{corollary}
\begin{remark}\label{r3} It is obvious that  the testing procedure attaining the lower bound in (\ref{49b})
is optimal among all truncated  testing procedures with
$\psi\in\Delta^N$. But  it only makes practical sense if
$$l_0=\min_{1\leq j\leq k}\sum_{i\not= j}\lambda_{ij}> 1 +R_0^N.$$

The reason is  that $l_0$ can be considered as "the $L(\chi,\psi)$"
function for a trivial sequential testing procedure
$(\chi_0,\psi_0,\phi_0)$ which, without taking any observations,
applies any decision rule $\phi_0$ such that $\phi _{0j}\leq
I_{\{\sum_{i\not = j}\lambda_{ij}=l_0\}}$ for any $j=1,\dots, k$. In
this case there are no observations ($N(\theta;\psi_0)=0$), $\chi_0$
is nothing, and it is easily seen that
$$L(\chi_0,\psi_0,\phi_0)=\sum_{j=1}^k\sum_{i\not =j}\lambda_{ij}\phi_{0j}=l_0.$$ Thus, the inequality
$$l_0\leq 1+R_0^N $$
means that the trivial testing procedure $(\chi_0,\psi_0,\phi_0)$ is not worse than the
best testing procedure with $\psi$ from $\Delta^N$.

Because of this, we may think that
$$V_0^N=\min\{l_0,1+R_0^N\}$$ is the  minimum value of $L(\chi,\psi)$ when taking no observations is permitted. It is obvious that this
is a particular case of (\ref{48}) with $n=1$, if we define
 $f_\theta^0\equiv 1$.
\end{remark}

\subsection{\hs General Stopping Rules}

In this section we characterize the structure of general sequential
testing procedures minimizing $L(\chi,\psi)$.

 Let us define for any stopping rule $\psi$ and any control policy
 $\chi$
\begin{equation}\label{50}
L_N(\chi,\psi)=\sum_{n=1}^{N-1}\int
s_n^{\psi,\chi}(nf_{\theta_1}^{n,\chi}+l_n^\chi)d\mu^n+\int c_N^{\psi,\chi}\left(Nf_{\theta_1}^{N,\chi}+l_N^\chi\right)d\mu^{N}.
\end{equation}
 This is the Lagrange-multiplier function
corresponding to  $\psi$ truncated  at $N$, i.e. the rule with the
components $\psi^N=(\psi_1,\psi_2,\dots,\psi_{N-1},1,\dots)$,
$L_N(\chi,\psi)=L(\chi,\psi^N)$.

Since $\psi^N$ is truncated, the results of the preceding section
apply, in particular, the inequalities of Theorem \ref{t3}.

The idea of what follows is to make $N\to\infty$, to obtain some
lower bounds for $L(\chi,\psi)$ from (\ref{46a}) - (\ref{46b}).
Obviously, we need that  $L_N(\chi,\psi)\to L(\chi,\psi)$ as
$N\to\infty$. A manner to guarantee this is using the following definition.

Let us denote by $\mathscr F$ the set of all strategies
($\chi,\psi$) such that
\begin{equation}\label{50a}
\lim_{n\to\infty} E_{\theta_i}^\chi
(1-\psi_1)\dots(1-\psi_{n})=0\quad\mbox{for any}\quad i=1,2,\dots,k.
\end{equation}
 It is easy to see that (\ref{50a}) is equivalent to
$$
P_{\theta_i}^\chi(\tau_\psi<\infty)=1\quad\mbox{for any}\quad
i=1,2,\dots,k
$$
(see (\ref{2a})).
\begin{lemma}\label{l3} For any  strategy $(\chi,\psi)\in \mathscr F$
$$
\lim_{N\to\infty}L_N(\chi,\psi)=L(\chi,\psi).
$$
\end{lemma}

\begin{proof} Practically coincides with that of Lemma 5.1 in \cite{NovikovIMF} (with $f_{\theta_1}^n$ instead of $f_{\theta_0}^n$), except that in order to show the convergence
$$
\int c_N^{\psi,\chi}l_N^\chi d\mu^N\to 0,\quad N\to\infty,
$$
we use the following estimate:
\begin{equation}\label{51}\int c_N^{\psi,\chi}l_N^\chi d\mu^N\leq
\max_{i\not =j}\lambda_{ij}\sum_{i=1}^k\int
c_N^{\psi,\chi}f_{\theta_i}^{N,\chi}d\mu^N=\max_{i\not =j}\lambda_{ij}\sum_{i=1}^kE_{\theta_i}^\chi c_N^\psi\to 0\end{equation}
as $N\to\infty$, because of (\ref{50a}).
\end{proof}

The second fact we need is about the behaviour of  the functions $V_r^N$
which participate in the inequalities of  Theorem \ref{t3}, as
$N\to\infty$.

\begin{lemma}\label{l4} For any $n\geq 1$ and for any  $N\geq n$
\begin{equation}\label{52a}
V_n^N\geq V_n^{N+1}.
\end{equation}
\end{lemma}

\begin{proof} Completely analogous to the proof of Lemma 5.2 \cite{NovikovIMF} (with $f_{\theta_1}^n$ instead of $f_{\theta_0}^n$).
\end{proof}

It follows from Lemma \ref{l4} that for any fixed $n\geq 1$ the
sequence $V_n^N$ is non-increasing. So, there exists
\begin{equation}\label{54}V_n= \lim_{N\to\infty}
V_n^N.
\end{equation}

Now,  passing to the limit, as
$N\to\infty$, in (\ref{46a}) and (\ref{46b}) with $\psi=\psi^N$, we have the following Theorem. The left-hand side of (\ref{46a}) tends to $L(\chi,\psi)$ by Lemma \ref{l3}. Passing to the limit on the right hand side of  (\ref{46a}) and in (\ref{46b}) is possible by Lebesgue's monotone convergence theorem, by virtue of Lemma \ref{l4}.

\begin{theorem}\label{t4} Let $(\chi,\psi)\in\mathscr F$ be any control-stopping strategy. Then for
any $r\geq 1$ the following inequalities hold
\begin{equation}\label{55}
L(\chi,\psi)\geq\sum_{n=1}^{r}\int
s_n^{\psi,\chi}(nf_{\theta_1}^{n,\chi}+l_n^\chi)d\mu^n+\int c_{r+1}^{\psi,\chi}\left((r+1)f_{\theta_1}^{r+1,\chi}+V_{r+1}^\chi\right)d\mu^{r+1}
\end{equation}
\begin{equation}\label{56}
\geq \sum_{n=1}^{r-1}\int
s_n^{\psi,\chi}(nf_{\theta_1}^{n,\chi}+l_n^\chi)d\mu^n
+\int
c_n^{\psi,\chi}\left(rf_{\theta_1}^{r,\chi}+V_{r}^\chi\right)d\mu^r,
\end{equation}
where
\begin{equation}\label{56a}
V_r=\min\{l_r,f_{\theta_1}^r+R_r\},
\end{equation}
being
\begin{equation}\label{56b}
R_r=R_r(x^{(r)},y^{(r)})=\min_{x_{r+1}}\int V_{r+1}(x^{(r+1)},y^{(r+1)})d\mu(y_{r+1}).
\end{equation}

In particular, for $r=1$, the following lower bound holds true:
\begin{equation}\label{57a}
L(\chi,\psi)\geq 1+\int V_1^\chi d\mu(y_1)\geq 1+R_0,
\end{equation}
where, by definition,
$$
R_0=\min_{x_1}\int V_1(x_1,y_1) d\mu(y_1).
$$
 \end{theorem}

Exactly as in \cite{NovikovIMF} (see Lemma 5.4 \cite{NovikovIMF}) it can be proved  that the right-hand side of (\ref{57a}) coincides
with
$$\inf_{(\chi,\psi)\in\mathscr F}L(\chi,\psi).$$
In fact, this is true for any $\mathscr F$ such that $(\chi,\psi)\in \mathscr F$ implies $L_N(\chi,\psi)\to L(\chi,\psi)$ as $N\to\infty$.

 The following theorem characterizes the structure of
optimal strategies.
\begin{theorem}\label{t5}  If there is a strategy $(\chi,\psi)\in \mathscr F$ such that
\begin{equation}\label{60}
L(\chi,\psi)=\inf_{(\chi^\prime,\psi^\prime)\in\mathscr F} L(\chi^\prime,\psi^\prime),
\end{equation}
then
\begin{equation}\label{60b}
I_{\{l_{r}^\chi<
f_{\theta_1}^{r,\chi}+R_{r}^{\chi}\}}\leq\psi_{r}^\chi\leq
I_{\{l_{r}^\chi\leq f_{\theta_1}^{r,\chi}+R_{r}^{\chi}\}}
\end{equation}
$\mu^r$-almost everywhere on $C_r^{\psi,\chi}$, and
\begin{equation}\label{60a}
\int V_{r+1}^{\chi}(y^{(r+1)})d\mu(y_{r+1})=R_{r}^{\chi}
\end{equation}
$\mu^r$-almost everywhere on $\bar C_r^{\psi,\chi}$,
for any $r=1, 2\dots$, where $\chi_1$ is defined in such a way that
\begin{equation}\label{60c}
\int V_1^{\chi} d\mu(y_1)=R_{0}.
\end{equation}

On the other hand, if a strategy $(\psi,\chi)$ satisfies (\ref{60b}) $\mu^r$-almost everywhere on $C_r^{\psi,\chi}$,
and satisfies (\ref{60a}) $\mu^r$-almost everywhere on $\bar C_r^{\psi,\chi}$, for any $r=1, 2\dots$, where $\chi_1$ is such that (\ref{60c}) is fulfilled,  and $(\psi,\chi)\in \mathscr F$, then (\ref{60}) holds.
\end{theorem}
\begin{proof}
Almost literally coincides with the proof of Theorem 5.5 \cite{NovikovIMF} (substituting $f_{\theta_0}^n$ by $f_{\theta_1}^n$), with the omission of the proof that $(\psi,\chi)\in \mathscr F$ in the "if"-part (see (76) and (77) in \cite{NovikovIMF}), because now it is a condition of Theorem \ref{t5}.

\end{proof}
\begin{remark}\label{r9}
Theorem \ref{t5} treats the optimality among strategies which take
at least one observation. If we allow to take no observations, there
is a possibility that the trivial testing procedure (see Remark
\ref{r3}) gives a better result. It is easy to see that this happens
if and only if
$$
l_0<1+R_0.
$$
\end{remark}

\section{\hs LIKELIHOOD RATIO STRUCTURE OF OPTIMAL STRATEGY}\label{s4}

In this section, we will give to the optimal strategy in Theorem
\ref{t5} an equivalent form related to the likelihood ratio process,
supposing that all the distributions given by $f_{\theta_i}$ are
absolutely continuous with respect to that given by $f_{\theta_1}$.
More precisely, we will suppose that for any $x$
\begin{equation}\label{6.1}
    \{y: f_{\theta_1}(y|x)=0\}\subset \bigcap_{i>1} \{y:
    f_{\theta_i}(y|x)=0\}.
\end{equation}

Let us start with defining the likelihood ratios:
$$
Z_n^r=Z_n^r(x^{(n)},y^{(n)})=\prod_{i=1}^n\frac{f_{\theta_r}(y_i|x_i)}{f_{\theta_1}(y_i|x_i)},\quad
r>1,
$$ and let $Z_n=(Z_n^2,\dots, Z_n^k)$.

 Let us introduce then the following sequence of
functions $\rho_r=\rho_r(z)$, $r=0,1,\dots$, where $z=(z_2,\dots
z_k)$.

Let
\begin{equation}\label{70a}
\rho_0(z)=g(z)\equiv \min_{j}\sum_{i\not =j }\lambda_{ij}z_i,
\end{equation}
where, by definition, $z_1\equiv 1$. Let for $r=1, 2, 3, \dots$,
recursively,
\begin{equation}\label{70}
\rho_r(z)=\min\left\{g(z),1+\min_{x}\int
f_{\theta_1}(y|x)\rho_{r-1}\left(z_2\frac{f_{\theta_2}(y|x)}{f_{\theta_1}(y|x)},\dots,
z_k\frac{f_{\theta_k}(y|x)}{f_{\theta_1}(y|x)}\right)d\mu(y)\right\}
\end{equation}
(we are supposing that all $\rho_r$, $r=0,1,2,\dots$ are
well-defined and measurable functions of $z$).
 It is easy to see that (see (\ref{48}),  (\ref{48a}))
$$
V_N^N=f_{\theta_1}^N\rho_0(Z_N),
$$
and for $r=N-1,N-2,\dots, 1$
\begin{equation}\label{72}
V_r^N=f_{\theta_1}^r\rho_{N-r}(Z_r).
\end{equation}
It is not difficult to see (very much like in Lemma \ref{l4}) that
\begin{equation*}
\rho_r(z)\geq \rho_{r+1}(z)
\end{equation*}
for any $r=0,1,2,\dots$, so there exists
\begin{equation}\label{72a}
\rho(z)=\lim_{n\to\infty}\rho_n(z).
\end{equation}
Using arguments similar to those used for obtaining Theorem
\ref{t4}, it can be shown, starting from (\ref{70}), that
\begin{equation}\label{71}
\rho(z)=\min\left\{g(z),1+R(z)\right\},
\end{equation}
where
\begin{equation}\label{72b}
R(z)=\min_{x}\int
f_{\theta_1}(y|x)\rho\left(z_2\frac{f_{\theta_2}(y|x)}{f_{\theta_1}(y|x)},\dots,z_k\frac{f_{\theta_k}(y|x)}{f_{\theta_1}(y|x)}\right)d\mu(y).
\end{equation}
Let us pass now to the limit, as $N\to\infty$, in (\ref{72}). We see
that
$$
V_k=f_{\theta_1}^k\rho(Z_k).
$$
Using this expression in Theorem \ref{t5} we get
\begin{theorem}\label{t6}
If  there exists a strategy $(\chi,\psi)\in \mathscr F$ such that
\begin{equation}\label{160}
L(\chi,\psi)=\inf_{(\chi^\prime,\psi^\prime)\in\mathscr F} L(\chi^\prime,\psi^\prime),
\end{equation}
then
\begin{equation}\label{160b}
I_{\{g(Z_r^\chi)< 1+R(Z_r^\chi)\}}\leq\psi_{r}^\chi\leq
I_{\{g(Z_r^\chi)\leq 1+R(Z_r^\chi)\}}
\end{equation}
$P_{\theta_0}^\chi$-almost sure on \begin{equation}\label{160bb}\{y^{(r)}:\; (1-\psi_1^\chi(y^{(1)}))\dots(1-\psi_{r-1}^\chi(y^{(r-1)}))>0\},\end{equation}
and
\begin{equation}\label{160a}
\int
f_{\theta_1}(y|\chi_{r+1})\rho\left(Z_r^{2,\chi}
\frac{f_{\theta_2}(y|\chi_{r+1})}{f_{\theta_1}(y|\chi_{r+1})},\dots,Z_r^{k,\chi}\frac{f_{\theta_k}(y|\chi_{r+1})}{f_{\theta_1}(y|\chi_{r+1})}\right)d\mu(y)=R(Z_r^{\chi})
\end{equation}
$P_{\theta_0}^\chi$-almost sure on \begin{equation}\label{160aa}\{y^{(r)}:\; (1-\psi_1^\chi(y^{(1)}))\dots(1-\psi_{r}^\chi(y^{(r)}))>0\},\end{equation} where $\chi_1$ is defined in such a way that
\begin{equation}\label{160c}
\int f_{\theta_1}(y|\chi_1)\rho\left(\frac{f_{\theta_2}(y|\chi_1)}{f_{\theta_1}(y|\chi_1)},\dots,\frac{f_{\theta_k}(y|\chi_1)}{f_{\theta_1}(y|\chi_1)}
 \right)d\mu(y)=R(1).
\end{equation}
On the other hand, if $(\chi,\psi)$ satisfies (\ref{160b}) $P_{\theta_0}^\chi$-almost sure on (\ref{160bb}) and
satisfies (\ref{160a}) $P_{\theta_0}^\chi$-almost sure on (\ref{160aa}), for any $r=1,2,\dots$, where $\chi_1$
satisfies (\ref{160c}),  and $(\chi,\psi)\in \mathscr F$, then $(\chi, \psi)$
satisfies (\ref{160}).
\end{theorem}

\section{\hs APPLICATION TO THE CONDITIONAL PROBLEMS}\label{s5}

In this section, we apply the results obtained in the preceding
sections to minimizing the average sample size
$N(\chi,\psi)=E_{\theta_1}^\chi \tau_\psi$ over all sequential
testing procedures  with error probabilities not exceeding some
prescribed levels (see Problems I and II in Section \ref{s1}).

Combining Theorems \ref{t1}, \ref{t2} and  \ref{t5}, we immediately
have the following solution to Problem I.
\begin{theorem}\label{t8} Let $(\chi,\psi)\in \mathscr F$ satisfy
the conditions of Theorem \ref{t5} with $\lambda_{ij}>0$, $i,j=1,\dots,k$,
$i\not =j$ (recall that $l_n$, $V_n$, and $R_n$ are functions of $\lambda_{ij}$), and let $\phi$ be any decision rule satisfying
(\ref{7aa}).

Then for any sequential testing procedure
$(\chi^\prime,\psi^\prime,\phi^\prime)\in \mathscr F$ such that
\begin{equation}\label{7.2}
    \alpha_{ij}(\chi^\prime,\psi^\prime,\phi^\prime)\leq
    \alpha_{ij}(\chi,\psi,\phi)\quad \mbox{for any}\quad
    i,j=1,\dots, k,\; i\not=j,
\end{equation}
it holds
\begin{equation}\label{7.3}
    N(\chi^\prime,\psi^\prime)\geq N(\chi,\psi).
\end{equation}

The inequality in (\ref{7.3}) is strict if at least one of the
inequalities in (\ref{7.2}) is strict.

If there are equalities in all of the  inequalities in (\ref{7.2})
and (\ref{7.3}), then $(\chi^\prime,\psi^\prime)$ satisfies
the condition of Theorem \ref{t5} as well (with $\chi^\prime$ instead of
$\chi$ and $\psi^\prime$ instead of   $\psi$).
\end{theorem}
\begin{proof} The only thing  to be proved is the last
assertion.

Let us suppose that
$$\alpha_{ij}(\chi^\prime,\psi^\prime,\phi^\prime)=
    \alpha_{ij}(\chi,\psi,\phi), \quad \mbox{for any}\quad
    i,j=1,\dots, k,\; i\not=j,$$
     and $$N(\chi^\prime,\psi^\prime)=
    N(\chi,\psi).$$
Then, obviously,
\begin{equation}\label{7.4}
    L(\chi,\psi,\phi)=L(\chi,\psi)=L(\chi^\prime,\psi^\prime,\phi^\prime)\geq
    L(\chi^\prime,\psi^\prime)
\end{equation}
(see (\ref{4})) and Remark \ref{r2}.

By Theorem \ref{t5}, there can not be strict inequality in the last
inequality in (\ref{7.4}), so
$L(\chi,\psi)=L(\chi^\prime,\psi^\prime)$. From Theorem \ref{t5} it
follows now that $(\chi^\prime,\psi^\prime)$ satisfies the conditions of Theorem \ref{t5} as well.
\end{proof}

Analogously, combining Theorems \ref{t1a}, \ref{t2} and  \ref{t5},
we also have the following solution to Problem II.
\begin{theorem}\label{t9} Let $(\chi,\psi)\in \mathscr F$ satisfy the conditions of Theorem \ref{t5} with $\lambda_{ij}=\lambda_i>0$ for any
$i=1,\dots k$ and for any $j=1,\dots, k$, and let $\phi$ be any
decision rule such that $$ \phi_{nj}\leq I_{\left\{\sum_{i\not
=j}\lambda_{i}f_{\theta_i}^n=\min_j\sum_{i\not
=j}\lambda_{i}f_{\theta_i}^n\right\}}$$ for any $j=1,\dots, k$ and
for any $n=1,2,\dots$.

Then for any sequential testing procedure
$(\chi^\prime,\psi^\prime,\phi^\prime)\in \mathscr F$ such that
\begin{equation}\label{7.102}
    \beta_{i}(\chi^\prime,\psi^\prime,\phi^\prime)\leq
    \beta_{i}(\chi,\psi,\phi)\quad \mbox{for any}\quad
    i=1,\dots, k,
\end{equation}
it holds
\begin{equation}\label{7.103}
    N(\chi^\prime,\psi^\prime)\geq N(\chi,\psi).
\end{equation}

The inequality in (\ref{7.103}) is strict if at least one of the
inequalities in (\ref{7.102}) is strict.

If there are equalities in all of the  inequalities in (\ref{7.102})
and (\ref{7.103}), then $(\chi^\prime,\psi^\prime)$ satisfies
the conditions of Theorem \ref{t5} with $\lambda_{ij}=\lambda_i$, $i,j=1,\dots, k$, $i\not=j$, as well (with $\chi^\prime$ instead of
$\chi$ and $\psi^\prime$ instead of   $\psi$).
\end{theorem}

\section{\hs ADDITIONAL RESULTS, EXAMPLES AND DISCUSSION}\label{s6}
\subsection{\hs Some general remarks}
\begin{remark}\label{r6}
The class $\mathscr F$ defined by (\ref{50a}) can be extended in such a way that Theorem \ref{t5} remains valid.
It can be defined as the class of all the strategies $(\chi,\psi)$ for which
\begin{equation}\label{61}
\lim_{n\to\infty} E_{\theta_i}^\chi
(1-\psi_1)\dots(1-\psi_{n})=0
\end{equation}
for at least $k-1$ different values of $\theta_i$. To see this it is sufficient to notice that for any strategy in this extended class
$$
L_N(\chi,\psi)\to L(\chi,\psi),\quad\mbox{as}\quad N\to\infty,
$$
because (see the proof of Lemma \ref{l3})
$$
\int c_N^{\psi,\chi}l_N^\chi d\mu^N\leq
\sum_{1\leq i\leq k,i\not=j}\lambda_{ij}\int
c_N^{\psi,\chi}f_{\theta_i}^{N,\chi}d\mu^N=\sum_{1\leq i\leq k,i\not=j}\lambda_{ij}E_{\theta_i}^\chi c_N^\psi\to 0,\quad N\to\infty,
$$
if $j$ corresponds to $\theta_j$ for which (\ref{61}) does not hold.

Obviously, Theorem \ref{t5} remains valid with this extension of $\mathscr F$.

Moreover, in the same way, Theorem \ref{t5} remains valid if $\mathscr F$ is defined as the class of all strategies $(\chi,\psi)$ for which
$$
L_N(\chi,\psi)\to\L(\chi,\psi),\quad N\to\infty.
$$
But the statistical meaning of this class is not clear, so we prefer for $\mathscr F$ one of the definitions above.
\end{remark}
\begin{remark}\label{r4}
In the same way as in the preceding sections, a more general problem than just minimizing $N(\theta_1;\chi,\psi)$  can
be treated (see
(\ref{8aa}) and Problems I and II thereafter).

Namely, we can minimize any convex combination of the average sample
numbers, or
$$
\sum_{i=1}^k c_iN(\theta_i;\chi,\psi),
$$
where $c_i\geq 0$, $i=1,\dots, k$, are arbitrary but fixed
constants. More exactly, if we modify  the definition of the
functions $V_r^N$ in (\ref{48}) to
\begin{equation}\label{48c}
V_{r-1}^N=\min\{l_{r-1},\sum_{i=1}^k
c_if_{\theta_i}^{r-1}+R_{r-1}^N\},
\end{equation}
for $r=N,\dots, 2$, being, as before,
$$
V_r=\lim_{N\to\infty}V_r^N,
$$
and, respectively, change (\ref{60b}) in Theorem \ref{t5} to
\begin{equation}\label{60bb}
I_{\{l_{r}^\chi<\sum_{i=1}^k
c_if_{\theta_i}^{r,\chi}+R_{r}^{\chi}\}}\leq\psi_{r}^\chi\leq
I_{\{l_{r}^\chi\leq \sum_{i=1}^k
c_if_{\theta_i}^{r,\chi}+R_{r}^{\chi}\}}
\end{equation}
then Theorem \ref{t5} remains valid. Theorems \ref{t3}, \ref{t6}, \ref{t8} and \ref{t9} can be modified respectively.
\end{remark}

\subsection{An example}
In this Section we show how our results can be applied to a concrete statistical model.

Let us suppose that any stage of our experiment is a regression experiment with a normal response. More specifically, we are supposing that the distribution of $Y$, given a value of the control variable $X$, is normal with mean value $\theta X$ and a know variance $\sigma^2$, say $\sigma^2=1$.

Thus,
\begin{equation}\label{90}
f_\theta(y|x)=\frac{1}{\sqrt{2\pi}}\exp\left\{-\frac{(y-\theta x)^2}{2}\right\},\quad -\infty<y<\infty
\end{equation}

For simplicity, let us take $k=2$ simple hypotheses, for example, $H_1:\,\theta=1$ and $H_2:\,\theta=2$, and suppose that the control variable takes only two values, say, $x=1$ and $x=2$.

Condition (\ref{6.1}) is fulfilled in an obvious way.

Let $\lambda_{12}>0$ and $\lambda_{21}>0$ two arbitrary constants.
We start defining
$$
\rho_0(z)=g(z)\equiv\min\{  \lambda_{12},\lambda_{21}z\},
$$
(see (\ref{70a})).

Next, we calculate
$$
\frac{f_{2}(y|x)}{f_{1}(y|x)}=\exp\{xy-3x^2/2\},
$$
and
$$
\rho_{n+1}(z)=\min\{g(z),1+\min_{x=1,2}\int_{-\infty}^\infty \rho_{n}(z\exp\{xy-3x^2/2\})\frac{\exp\{-(y-x)^2/2\}}{\sqrt{2\pi}}dy,\quad
$$
for $n=0,1,2,\dots$ (see (\ref{70})).

Let $\rho(z)=\lim_{n\to\infty}\rho_n(z)$, and $$R(z)=\min_{x=1,2}\int_{-\infty}^\infty \rho(z\exp\{xy-3x^2/2\})\frac{\exp\{-(y-x)^2/2\}}{\sqrt{2\pi}}dy.$$

Now, by Theorem \ref{t6}, an optimal strategy will be defined on the basis of the likelihood ratio process
$$
Z_n=\exp\{\sum_{i=1}^n (X_i Y_i-3X_i^2/2)\},
$$
being the optimal stopping time $\tau=\min\{n:g(Z_n)\leq 1+R(Z_n)\}$, whereas at each stage $n=1,2,\dots$
the next control value $X_{n+1}=x$ ($x=1$ or $x=2$) is defined in such a way that
$$
R(Z_n)=\int_{-\infty}^\infty \rho(Z_n\exp\{xy-3x^2/2\})\frac{\exp\{-(y-x)^2/2\}}{\sqrt{2\pi}}dy,
$$
starting from $X_1$ defined as $x$ ($x=1$ or $x=2$) for which
$$
R(1)=\int_{-\infty}^\infty \rho(\exp\{xy-3x^2/2\})\frac{\exp\{-(y-x)^2/2\}}{\sqrt{2\pi}}dy.
$$
When the test terminates at some stage $\tau=n$, we should reject $H_1$, if  $\lambda_{21}Z_n\geq \lambda_{12}$, and accept $H_1$ otherwise (see Theorem \ref{t3}).

One can vary the error probability levels of this test by changing the values of $\lambda_{12}$ and $\lambda_{21}$.

\subsection{\hs Bayesian testing of multiple hypotheses}\label{ss6.2}
In this section we characterize the structure of Bayesian multiple hypothesis tests.

Let $\pi_i>0$, $i=1,\dots, k$ be prior probabilities of $H_i$, $i=1,\dots, k$, respectively, $\sum_{i=1}^k\pi_i=1$, and let $w_{ij}\geq 0$, $i,j=1,\dots, k$,  be some losses due to incorrect decisions (we assume that $w_{ii}=0$ for any $i=1,\dots, k$).
Then, for any sequential testing procedure $(\chi,\psi,\phi)$, we define the Bayes risk as
\begin{equation}\label{62}
r(\chi,\psi,\phi)=\sum_{i=1}^k \pi_i\left(cE_{\theta_i}^\chi \tau_\psi +\sum_{j=1}^kw_{ij}\alpha_{ij}(\chi,\psi,\phi)\right),
\end{equation}
where $c>0$ is some unitary observation cost (cf. Section 9.4 of \cite{Zacks}, see also Chapter 5 of \cite{Ghosh} for a more general sequential Bayesian decision theory, both monographs treating non-controlled experiments). Let us call {\em Bayesian} any testing procedure $(\chi,\psi,\phi)$ minimizing (\ref{62}).

In this section, we show that the Bayesian testing procedures always exist, and characterize the structure of both truncated and non-truncated Bayesian testing procedures for the controlled experiments.

To formulate our results, we use the notation of Sections \ref{s1} - \ref{s5}, but we have to re-define some elements have been defined therein.

First of all, it is easy to see from Theorem \ref{t2} that the optimal decision rule $\phi$ has the following form. Let
\begin{equation}\label{63}
l_n=\min_{1\leq j\leq k}\sum_{i=1}^k\pi_iw_{ij}f_{\theta_i}^n.
\end{equation}
(cf. (\ref{6a})). Then the decision rule $\phi$ is optimal ($\inf_{\phi^\prime}r(\chi,\psi,\phi^\prime)=r(\chi,\psi,\phi^\prime)$ for any $\chi$ and $\psi$) if
\begin{equation}\label{64}
\displaystyle
\phi_{nj}\leq I_{\{\sum_{i=1}^k\pi_iw_{ij}f_{\theta_i}^n=l_n\}}
\end{equation}
for any $j=1,\dots, k$ and for any $n=1,2,\dots$ (see Theorem \ref{t2}).

Let $\Pi$ be the prior distribution defined by $\pi_i$, $i=1,\dots,k$, and let, by definition,
$$
f_\Pi^{n}=\sum_{i=1}^k \pi_if_{\theta_i}^n
$$
for any $n=1,2,\dots$.

For any $N=1,2,\dots$ let us define
\begin{equation}\label{65}
V_N^N=l_N,
\end{equation}
and for any $n=N-1,N-2,\dots, 1$, recursively,
\begin{equation}\label{66}
    V_n^N=\min\{l_n, cf_\Pi^n+R_n\},
\end{equation}
where
\begin{equation}\label{67}
    R_n^N=R_n^N(x^{(n)},y^{(n)})=\min_{x_{n+1}}\int V_{n+1}^N(x_1,\dots, x_{n+1};y_1,\dots,y_{n+1})d\mu(y_{n+1}).
\end{equation}
Let also
\begin{equation}\label{68}
    R_0^N=\min_{x_1}\int V_1^N(x_1;y_1)d\mu(y_1).
\end{equation}
The following Theorem characterizes Bayesian procedures with truncated stopping rules and can be  proved in exactly the same way as Corollary \ref{c1}.
\begin{theorem}\label{t10} Let $\chi$ be any control policy, $\psi\in\Delta^N$ be any (truncated) stopping rule and $\phi$ any decision rule satisfying (\ref{64}) for any $j=1,\dots, k$ and for any $n=1,2,\dots$. Then
\begin{equation}\label{69}
          r(\chi,\psi,\phi)\geq c+R_0^N.
\end{equation}
There is an equality in (\ref{69}) if and only if
\begin{equation}\label{80}
I_{\{l_{n}^\chi< cf_\Pi^{n,\chi}+R_n^{N,\chi} \}}\leq\psi_{n}^\chi\leq I_{\{l_{n}^\chi\leq cf_\Pi^{n,\chi}+R_n^{N,\chi} \}}
\end{equation}
$\mu^n$-almost everywhere on $C_n^{\psi,\chi}$
and
\begin{equation}\label{81}
R_n^{N,\chi}(y^{(n)})=\int V_{n+1}^{N,\chi}(y^{(n+1)})d\mu(y_{n+1})
\end{equation}
 $\mu^n$-almost everywhere on $\bar C_{n}^{\psi,\chi}$, for any $n=1,\dots, N-1$, and, additionally,
 \begin{equation}\label{82}
R_0^N=\int V_{1}^{N}(\chi_1;y_1)d\mu(y_{1}).
\end{equation}
\end{theorem}

Let now $V_n=\lim_{N\to\infty} V_n^N$, $n=1,2,\dots$. Respectively, $R_n=\lim_{N\to\infty}R_n^N$, $n=0,1,2,\dots$.
\begin{theorem}\label{t11} Let $\chi$ be any control policy, $\psi$ any  stopping rule, and $\phi$ any decision rule satisfying (\ref{64})  for any $j=1,\dots, k$ and for any $n=1,2,\dots$. Then
\begin{equation}\label{69a}
           r(\chi,\psi,\phi)\geq c+R_0.
\end{equation}
There is an equality in (\ref{69a}) if and only if
\begin{equation}\label{80a}
I_\{l_{n}^\chi< cf_\Pi^{n,\chi}+R_n^{\chi} \}\leq\psi_{n}^\chi\leq I_{\{l_{n}^\chi\leq cf_\Pi^{n,\chi}+R_n^{\chi}\}}
\end{equation}
$\mu^n$-almost everywhere on $C_n^{\psi,\chi}$
and
\begin{equation}\label{81a}
R_n^{\chi}(y^{(n)})=\int V_{n+1}^{\chi}(y^{(n+1)})d\mu(y_{n+1})
\end{equation}
 $\mu^n$-almost everywhere on $\bar C_{n}^{\psi,\chi}$, for any $n=1, 2\dots $, and, additionally,
 \begin{equation}\label{82a}
R_0=\int V_{1}(\chi_1;y_1)d\mu(y_{1}).
\end{equation}
\end{theorem}
\begin{proof}
First of all we need to prove that (\ref{69a}) holds for {\em any} strategy $(\chi,\psi)$.  Obviously, it suffices to prove this only for such $(\chi,\psi)$ that $r(\chi,\psi,\phi)<\infty$. But this latter fact implies, in particular, that $\sum_{i=1}^k \pi_i E_{\theta_i}^\chi \tau_\psi<\infty$ (see (\ref{62})). Because $\pi_i>0$ for any $i=1,\dots k$, it follows that $(\chi,\psi)$ satisfies (\ref{50a}), so
$$
r(\chi,\psi^N,\phi)\to r(\chi,\psi,\phi),\quad N\to\infty,
$$
where $\psi^N$, by definition, is $(\psi_1,\psi_2,\dots,\psi_{N-1},1,\dots)$ (see the proof of Lemma \ref{l3}).

The rest of the proof of the "only if"-part is completely analogous to the corresponding part of the proof of Theorem \ref{t5} (or Theorem 5.5 \cite{NovikovIMF}).

To prove the "if"-part, first it can be shown, analogously to the proof of Theorem 5.5 \cite{NovikovIMF}, that
\begin{equation}\label{83}
\sum_{n=1}^{r}\int
s_n^{\psi,\chi}(cnf_\Pi^{n,\chi}+l_n^\chi)d\mu^n+\int
c_{r+1}^{\psi,\chi}\left(c(r+1)f_\Pi^{r+1,\chi}+V_{r+1}^\chi \right)d\mu^{r+1}=c
+R_0,
\end{equation}
for any $r=0,1,2,\dots$, if $(\psi,\chi)$ satisfies (\ref{80a}) -- (\ref{82a}).

Because $c>0$, we have from (\ref{83}), in particular, that
$$
\sum_{i=1}^k\pi_iP_{\theta_i}^\chi(\tau_\psi\geq r+1)=\int
c_{r+1}^{\psi,\chi}f_\Pi^{r+1,\chi}d\mu^{r+1}\leq\frac{c+R_0}{c(r+1)}\to 0\quad\mbox{as}\quad r\to\infty.
$$
Because $\pi_i>0$ for all $i=1,\dots, k$, this implies that for $(\chi,\psi)$ (\ref{50a}) is fulfilled.
It follows from (\ref{83}) now that
$$
\lim_{r\to\infty}\sum_{n=1}^{r}\int
s_n^{\psi,\chi}(cnf_\Pi^{n,\chi}+l_n^\chi)d\mu^n=r(\chi,\psi,\phi)\leq c+R_0.
$$
Along with (\ref{69a}) this gives that $r(\chi,\psi,\phi)=c+R_0$, i.e. there is an equality in (\ref{69a}).
\end{proof}
\subsection{\hs Experiments without control}

In this section we draw consequences for statistical experiments without control.

Let us suppose that the density of $Y$ given $X$ does not depend on $X$: $f_\theta(y|x)\equiv f_\theta(y)$ for any $y$ and for any $\theta$, meaning that there is no way to control the flow of the experiment, and the observations $Y_1,Y_2,\dots$ are independent and identically distributed (i.i.d.) random "variables" with probability "density" function $f_\theta(y)$. We can incorporate this particular case in the above scheme of controlled experiments thinking that there is some (fictitious) unique value of control variable at each stage of the experiment, thus, being any control policy trivial.

  Because of this,   any (sequential) testing procedure has in effect only two components in this case: a stopping rule $\psi$ and a decision rule $\phi$. So we use the notation of section \ref{ss6.2}, simply omitting any mention of the control policy. For example, for any testing procedure $(\psi,\phi)$ the Bayesian risk (\ref{62}) is now:
\begin{equation}\label{62a}
r(\psi,\phi)=\sum_{i=1}^k \pi_i\left(cE_{\theta_i} \tau_\psi +\sum_{j=1}^kw_{ij}\alpha_{ij}(\psi,\phi)\right).
\end{equation}
Respectively,
$f_{\theta}^n=f_{\theta}^n(y^{(n)})=\prod_{i=1}^nf_{\theta}(y_i)$ in (\ref{63}) now, and the functions $V_n^N$, $R_n^N$, $V_n$, $R_n$, etc. of the preceding section are all functions of $y^{(n)}$ only.

Theorem \ref{t11} of section \ref{ss6.2} transforms now to
\begin{theorem}\label{t12} Let  $\psi$ be any stopping rule and $\phi$ any decision rule satisfying (\ref{64}) for any $j=1,\dots, k$ and for any $n=1,2,\dots$. Then
\begin{equation}\label{69b}
          r(\psi,\phi)\geq c+R_0.
\end{equation}
There is an equality in (\ref{69b}) if and only if
\begin{equation}\label{80b}
I_{\{l_{n}< cf_\Pi^{n}+R_n \}}\leq\psi_{n}^\chi\leq I_{\{l_{n}\leq cf_\Pi^{n}+R_n \}}
\end{equation}
$\mu^n$-almost everywhere on $C_n^{\psi}$ for any $n=1,2,\dots$,  where
$$
R_n=R_n(y_1,\dots, y_n)=\int V_{n+1}(y_1,\dots,y_{n+1})d\mu(y_{n+1}),
$$
being, for any $n=1,2,\dots$, $V_n(y^{(n)})=\lim_{N\to\infty}V_n^N(y^{(n)})$, where
$V_N^N\equiv l_N$, and
$$
V_n^N(y^{(n)})=\min\{l_n(y^{(n)}),cf_\Pi^n(y^{(n)})+\int V_{n+1}^N(y^{(n+1)})d\mu(y_{n+1)}\}
$$
for any $n=N-1,\dots,1$, $N=1,2,\dots$
\end{theorem}
In particular, this Theorem gives all solutions to the
problem of Bayesian testing of multiple simple hypotheses for independent and
identically distributed observations
when the cost of observations is linear (see Section 9.4 of
\cite{Zacks} and suppose  that $K(X_1,\dots,X_n)\equiv n$ therein).

In the particular case of two hypotheses ($k=2$) a Bayesian test of Theorem \ref{t12} given by
$$
\psi_{n}^\chi=I_{\{l_{n}\leq cf_\Pi^{n}+R_n \}},\quad n=1,2,\dots,
$$
has the form of the Sequential Probability Ratio Test (SPRT, see \cite{WaldWolfowitz}), being all other Bayesian tests (\ref{80b}) randomizations at its boudaries (see  \cite{NovikovIJPAM} for closely related results).

 \vspace{4mm}

\section*{ACKNOWLEDGEMENTS}
\small
The author thanks the anonymous referees for reading the manuscript very carefully, and for their valuable comments and suggestions.

The author greatly appreciates
the support of the Autonomous
 Metropolitan University, Mexico City, Mexico, where this work was
 done, and the support of the National System of Investigators (SNI)
 of
 CONACyT, Mexico.

This work is also partially supported
by Mexico's CONACyT Grant no.
CB-2005-C01-49854-F.\vspace{3mm}

\footnotesize
\begin{flushright}
(Received April 18, 2008.)\,\
\rule{0mm}{0mm}
\end{flushright}

\vspace*{2mm}

{\mi
\begin{flushright}
\begin{minipage}[]{124mm}
{Andrey Novikov, Departamento de
Matem\'aticas, Universidad Aut\'onoma
Metropolitana - Unidad Iztapalapa, San
Rafael Atlixco 186, col. Vicentina,
C.P. 09340, M\'exico D.F., M\'exico
\\ e-mail: {\tt an@xanum.uam.mx}
\\{\tt http://mat.izt.uam.mx/profs/anovikov/en}}
\end{minipage}
\end{flushright}
}
\end{document}